\documentclass[10pt]{amsart}
\usepackage{mystyle}

\title[Einsiedler-Polo-Kadyrov Inequality]{Pressure Inequalities for Gibbs Measures of Countable Markov Shifts}
\author{Ren\'e R\"uhr}
\address{Weizmann Institute
{\tt rene.ruhr@weizmann.ac.il}
}
\begin{document}

\begin{abstract}
	We provide a quantification of the uniqueness of Gibbs measure for topologically mixing countable Markov shifts with locally H\"older continuous potentials. Corollaries for speed of convergence for approximation by finite subsystems are also given.
\end{abstract}

\maketitle
\setcounter{tocdepth}{1}
% \tableofcontents

\section{Statement of results}
Consider a countable Markov shift $(T,X,A)$ where 
\[X=\{ x=(x_0,x_1,\dots)\in S^{\N_0} : t_{x_ix_j}=1\}\]
for a countable sets of states $S$, $A=(t_{ij})_{i,j\in S}$ denotes the transition matrix and $T:X\to X$ is the left shift, which we assume to be topologically mixing. Given a potential $\phi:X\to \R$, a $T$-invariant probability measure $m$ is called a Gibbs measure for $\phi$ if there exist constants $C,P$ such that for any cylinder $[a_0,\dots,a_{n-1}]=\{x\in X: x_i=a_i, 0\leq i<n\}$ and for all $x\in[a_0,\dots,a_{n-1}]$ it holds that
\[
\frac{1}{C}\leq \frac{m([a_0,\dots,a_{n-1}])}{e^{\phi_n(x)-nP}} \leq C
\]
where $\phi_n=\phi+\phi\circ T+\dots+\phi\circ T^{n-1}$.
Fix some $\theta\in(0,1)$ and introduce the metric $d(x,y)=\theta^{\min\{i: x_i\neq y_i\}}$ on $X$. For what follows we assume that $\phi$ is locally H\"older continuous, that is, there is a constant $A_\phi>0$ so that
\[
\var_n\phi\leq A_\phi \theta^n
\quad\text{ for all } n\geq1 \text{ where}\quad\var_n\phi=\sup_{x,y}\{\phi(x)-\phi(y):x_i=y_i, i<n\}.
\]
In particular, $\phi$ is H\"older summable, $\sum_{n\geq1}\var_n(\phi)<\infty$.
The Gurevich-Sarig pressure of $(T,X)$ with respect to $\phi$ is
\[
P_{GS}(\phi)=\lim_{n\to\infty}\frac{1}{n}\log Z_n(\phi,a)\quad \quad\quad\quad\quad
Z_n(\phi,a)=\sum_{ x\in \text{Fix}_n}e^{\phi_n(x)}\id_{[a]}(x)
\]
where $\text{Fix}_n=\{x\in X: T^nx=x\}$ and $a\in S$ arbitrary. The limit exists (possibly infinite) and is independent of $a$ (\cite{sarig1999thermodynamic}[Theorem~1]).
For an arbitrary $T$-invariant probability $\mu$ on $X$ introduce also the metric pressure
\[
P_\mu(\phi)=h_\mu(T)+\mu(\phi) 
\]
where $h_\mu(T)$ denotes the Kolmogorov-Sinai entropy of $T$ and $\mu(\phi)=\int \phi d\mu$. This will always be well defined if $-\mu(\phi)<\infty$.
 The variational principle \cite{sarig1999thermodynamic} in this context states that if $\|L_\phi 1\|_\infty<\infty$, where $L_\phi f(x)=\sum_{Ty=x} e^{\phi(y)}f(y)$ then 
\[P_{GS}(\phi)=\sup_\mu P_\mu(\phi)\]
 where the supremum runs over all $T$-invariant Borel probability measures $\mu$ for which $-\mu(\phi)<\infty$. 
 In our setting where there is a Gibbs measure, $P_{GS}(\phi)<\infty$ (see Remark~\ref{rem:finitepressure}) and we may remove the condition $-\mu(\phi)<\infty$ for the definition of $P_\mu(\phi)$ by defining instead $P_\mu(\phi)=\mu(I_\mu+\phi)$ where $I_\mu$ is the information function associated to $\mu$ (see Remark~\ref{rem:onesidedintegrable}). Theorem~\ref{thm:main_theorem} below will employ this convention.

For a continuous function $f:X\to\R$ let
\[
\|f\|_\theta=\sup\{|f(x)-f(y)|/d(x,y): x_0=y_0\}
\]
be the Lipschitz constant of $f$. Let $\cL$ denote the Banach space of functions $f$ for which $\|f\|_{\cL}=\|f\|_\theta+\|f\|_\infty$ is finite.

In this note, we establish the following.
\begin{theorem}\label{thm:main_theorem}
	Assume that $(T,X)$ is a topological mixing countable Markov shift and that $\phi:X\to\R$ is a locally H\"older continuous potential. Assume that $m$ is a Gibbs measure for $\phi$. Then there is a constant $a>0$ such that for any $T$-invariant probability $\mu$ and any $f\in\cL$,
	\begin{equation}
		|m(f)-\mu(f)|\leq a \|f\|_{\cL}\left(P_{GS}(\phi)-P_\mu(\phi)\right)^{\frac12}.
	\end{equation}
\end{theorem}
This theorem is a generalization of a theorem of Kadyrov \cite{kadyrov2015effective}, which proves Theorem~\ref{thm:main_theorem} when $S$ is finite and $\phi=0$. Note that in this setting, the expression in the square root in the right-hand side becomes $h_{\text{top}}(T)-h_\mu(T)$. 
\begin{remark}
	M.\ Einsiedler gave the original argument we employ for the $2x$ map on the torus. The argument has been worked out in detail and extended for the cat map in F.\ Polo's thesis \cite{polo} and is generalized to hold for subshifts of finite type by S. Kadyrov \cite{kadyrov2015effective},\cite{kadyrov2016effective}. A treatment for $p$-adic diagonalizable actions on a compact homogeneous space is given by the author \cite{ruhr2016effectivity}. Using a very different argument, a similar estimate is also obtained for $p$-adic actions on non-compact homogeneous spaces by I.\ Khayutin \cite{khayutin2017large}.

Einsiedler's argument is in fact finitary in the sense that $h_\mu(T)$ may be replaced by $\frac{1}{n}H_\mu(\xi_0^{n-1})$ where $\xi=\{[a]\}_{a\in S}$ is the partition of $X$ into principal cylinders, $\xi_0^{n-1}=\bigvee_{j=0}^{n-1} T^{-j}\xi$ and
$H_\mu(\eta)=-\sum_{P\in \eta}\mu(P)\log\mu(P)$ for a partition $\eta$.
 We also establish this variant generalizing \cite{kadyrov2016effective} which treats the subshift of finite type setting.
\end{remark}

\begin{theorem}\label{thm:finitary_theorem}
	Assume that $(T,X)$ is a topological mixing countable Markov shift and that  $\phi:X\to\R$ is a locally H\"older continuous potential. Assume that $m$ is a Gibbs measure for $\phi$. Then there exists a constant $b>0$ such that for any $T$-invariant probability $\mu$ for which $-\mu(\phi)<\infty$ and $H_\mu(\xi)<\infty$ and for any $f\in\cL$, $n\in\N_{>0}$,
	\begin{equation}
		|m(f)-\mu(f)|\leq b\|f\|_{\cL}\left(\theta^{n}+\left( 2\theta^{n}\sum_{k\geq n}\var_k(\phi)+ P_{GS}(\phi)-\left(\mu(\phi)+H_\mu(\xi|\xi_0^{n-1})\right)\right)^{\frac12}\right).
	\end{equation}
If for some $\ell\geq1$ it holds that $\var_\ell(\phi)=0$, then for any $n\geq3\ell$ also
\[
|m(f)-\mu(f)| \leq b\|f\|_{\cL} \left( \theta^{\lfloor\frac{n-\ell}{2}\rfloor} +  \left(  P_{GS}(\phi) - \left( \mu(\phi)+\frac{1}{n} H_\mu(\xi_0^{n-1})  \right) + \frac{\ell}{n-\ell}H_\mu(\xi)  \right)^{\frac12} \right).
\]

\end{theorem}

\begin{remark}\label{rem:finitepressure}
We emphasize that it is the insight of Sarig \cite{sarig2003existence} that countable Markov shifts of the above type that have a Gibbs measure enjoy the Big Images and Preimages property (BIP), that is, there exists a finite collection of states $S'\subset S$ such that for any $a\in S$ there is $b,c\in S'$ such that $t_{ba}t_{ac}=1$. BIP on the other hand implies that $(X,T)$ essentially behaves like a subshift of finite type, which allows us to extend the results of Kadyrov to this setting with little effort. 
Another result of \cite{sarig2003existence} is that the existence of a Gibbs measure implies that $P_{GS}(\phi)=\lim_{n\to\infty}\frac{1}{n}\log \sum_{ x\in \text{Fix}_n}e^{\phi_n(x)}$ and that $P_{GS}(\phi)$ is finite, so that the above theorems are non-vacuous. Finally, if there is a Gibbs measure, then it is unique, see Theorem~4.9 \cite{sarig2009lecture}.
\end{remark}
\begin{remark}\label{rem:onesidedintegrable}
	We also need to mentioned that the Gibbs measure $m$ is not assumed to be an equilibrium measure, that is, $h_m(T)+m(\phi)\neq P_{GS}(\phi)$ which happens if $-m(\phi)=h_m(T)=\infty$. In this case, one can still work with $m(I_m+\phi)$ $(=P_{GS}(\phi))$ where $I_m$ is the information function. In fact, $I_\mu+\phi$ will be one-sided integrable for any $T$-invariant probability measure $\mu$ by Corollary 2 in \cite{sarig1999thermodynamic}.
	The case in which $m$ is an equilibrium measure for a H\"older summable potential, $P_m(\phi)=P_{GS}(\phi)<\infty$ and $(X,T)$ has the BIP property then it is in fact a Gibbs measure (\cite{sarig2003existence},\cite{mauldin2001gibbs}). 
\end{remark}

\begin{acknowledgments}
Work on this research matter received partial support by ISF grant 1149/18.
\end{acknowledgments}

\section{Proof of Theorem~\ref{thm:main_theorem}}

The proof of Theorem~\ref{thm:main_theorem} follows essentially verbatim \cite{kadyrov2015effective} by using Sarig's Generalized Ruelle's Perron-Frobenius theorem (GRPF) \cite{sarig2001thermodynamic}.
As mentioned above, we have the following \cite{sarig2003existence}: $(X,T)$ has a Gibbs measure $m$ for a H\"older summable potential $\phi$ if and only if $(X,T)$ satisfies BIP and $P_{GS}(\phi)<\infty.$ BIP in particular implies that $\phi$ is positive recurrent so that the GPRF guarantees the existence of a positive continuous function $h:X\to\R$ and a conservative measure $\nu$ such that $L_\phi^*\nu=\lambda \nu$ and $L_\phi h = \lambda h$ with $\lambda=\exp P_{GS}(\phi)$ where $L_\phi$ is Ruelle's operator $L_\phi f(x)=\sum_{Ty=x} e^\phi(y)f(y)$. 
$h$ can be normalized such that the measure $dm_0=hd\mu$ is probability. If $m$ is a Gibbs measure then $m=m_0$ (see Theorem~ 4.9 \cite{sarig2009lecture}). BIP further implies that $h$ is uniformly bounded away from $0$ and infinity, Corollary~2 \cite{sarig2003existence}, in particular $\nu$ is finite.

By a theorem of Aaronson and Denker \cite{aaronson2001local} it follows that BIP implies that there exists $\kappa\in(0,1)$ and $c>0$ such that for any $f\in\cL$
\begin{equation}\label{eq:spectralgap}
	\|L_\phi^n f-h\nu(f)\|_{\cL}\leq c \kappa^n\|f\|_{\cL},
\end{equation}
see Theorem~5.8, \cite{sarig2009lecture}.

Let $\xi$ denote the partition $\{[a]\}_{a\in S}$ so that $\xi_0^\infty=\cB_X$. For a measure $\mu$ let $g_\mu=\frac{d\mu}{d\mu\circ T}$ where the measure $\mu\circ T$ is defined by $\mu\circ T(f)=\sum_{a\in S}\int_{T[a]}f(ax)d\mu(x)$. Then the information function $I_\mu(x)=I_\mu(\xi|\xi_1^\infty)(x)=-\log \mu^{\xi_1^\infty}_x([x_0])$ satisfies Ledrappier's formula 
\[
I_\mu=-\log g_\mu.
\] 
As noted, $\phi$ is positive recurrent so that by Theorem~6,\cite{sarig1999thermodynamic}
\begin{equation}\label{eq:cohom}
I_m=-[\phi+\log h-\log h\circ T-P_{GS}(\phi)].
\end{equation}
We further note that \eqref{eq:cohom} implies that $I_m$ and hence also $m_x^{\xi_1^\infty}([x_0])$ is defined everywhere (and not just $m$-a.e.). It also follows that $m(I_m+\phi)$ is defined and equal to $P_{GS}(\phi)$, and that $\mu(I_m)=P_{GS}(\phi)-\mu(\phi)$.

Consider the transfer operator $\widehat{T}f=\frac{dT_*m_f}{dm}=h^{-1}L_\phi(hf)$, $dm_f=fdm$ that satisfies
\(
\widehat{T}f\circ T=\E_m[f|\xi_1^\infty].
\)
We remark that this equality holds everywhere as both sides of the equation are continuous functions.
Let $f_n=\widehat{T}^nf=h^{-1}L_\phi^n(hf)$, which converges to $m(f)$ by the GRPF.
Note that $\E_m[f|\xi_1^\infty](x)=\sum_{s\in S}f(sTx)m^{\xi_1^\infty}_x([s])$ and similarly for $\mu$ where $sTx=(s,x_1,\dots)$.

Introduce the probability vectors $p_s(x)=m^{\xi_1^\infty}_x([s])$, $q_s(x)=\mu^{\xi_1^\infty}_x([s])$ and the Kullback-Leibler convergence of two probability vectors $p,q$
\[
D_p(q)=-\sum q_i\log\frac{p_i}{q_i}
% \]
\quad\text{ satisfying }\quad
% \[
\|q-p\|_1\leq \sqrt{2D_p(q)}
\]
(Pinsker's inequality).
Write $m(f)-\mu(f)=\mu(\lim f_n-f_0)$. Telescoping and interchanging summation and integration (Fubini, justified by inequality~\eqref{eq:spectralgap}) gives
\(
m(f)-\mu(f)=\sum_{n\geq0} \mu(f_{n+1}-f_n).
\)
We find,
\[
\mu(f_{n+1}-f_n)=\mu(f_{n+1}\circ T-f_n)=\mu(\E_m[f_n|\xi_1^\infty]-\E_\mu[f_n|\xi_1^\infty])
\]
\[
=\int d\mu(x)\sum_{s\in S}f_n(sTx)(p_s(x)-q_s(x))=\int d\mu(x)\sum_{s\in S}(f_n(sTx)-m(f))(p_s(x)-q_s(x))
\]
where we added zero in the last equality since $\sum_s(p_s-q_s)=0$. 
We bound $f_n-m(f)=h^{-1}(L^n_\phi(fh)-h\nu(fh))$ by its sup norm using that $h$ is bounded from below from zero. 
As the sup norm is dominated by $\|\cdot\|_\cL$, since $h$ is also bounded from above and finally $\|hf\|_\cL\leq\|h\|_\infty\|f\|_\cL$ we conclude from the spectral gap, inequality~\eqref{eq:spectralgap} that $|(f_n(sTx)-m(f))|\leq c\kappa^n\|h^{-1}\|_\infty\|h\|_\infty\|f\|_\cL$.

By Pinsker's inequality applied to the remaining sum gives
\begin{equation}\label{eq:boundKL}
|\mu(f_{n+1}-f_n)|\leq \|f_n-m(f)\|_\infty\mu(\|p-q\|_1)\leq c\kappa^n\|h^{-1}\|_\infty\|h\|_\infty\|f\|_\cL \left(2\mu(D_p(q))\right)^{\frac12}.
\end{equation}
The integral reduces to the difference of pressures,
\[
	\mu(D_{p(\cdot)}(q(\cdot)))=\mu(-\log p(\cdot)+\log q(\cdot))=\mu(I_m-I_\mu)=P_{GS}(\phi)-\mu(\phi+I_\mu)=P_{GS}(\phi)-P_\mu(\phi)
\]
where we have used properties of conditional measures in the first equality and equation~\eqref{eq:cohom} in the third equality.
The theorem now follows with $a=\frac{\sqrt{2}c}{1-\kappa}\|h^{-1}\|_\infty\|h\|_\infty$.

\section{Proof of Theorem~\ref{thm:finitary_theorem}}
We can obtain Theorem~\ref{thm:finitary_theorem} by modifying the proof of Theorem~\ref{thm:main_theorem} as follows.
We replace the probability vector $q_s(x)$ by $q^m_s(x)=\mu^{\xi_1^{m-1}}_x([s])$.
 Then
\begin{equation}\label{eq:approx_m}
\left| E_\mu(f_n|\xi_1^{m-1})(x)-\sum_{s\in S} f_n(sTx)q^m_s(x) \right|
=\left|\sum_{s\in S} \int_{[s]} (f_n(y)-f_n(sTx)) d\mu^{\xi_1^{m-1}}_x(y) \right|
\leq \theta^{m}\|f_n\|_\theta
\end{equation} 
 since $d(y,sTx)\leq \theta^{m-1}$ for any $y\in[s]\cap[x]_{\xi_1^{m-1}}$.
Hence 
\[
\left|\mu(\E_m[f_n|\xi_1^\infty]-\E_\mu[f_n|\xi_1^\infty])\right|
\leq \theta^{m}\|f_n\|_\theta + \left|\mu\left(\E_m[f_n|\xi_1^\infty]-\sum_{s\in S} f_n(sT\cdot)q^m_s(\cdot)\right)\right|
\]
and the proof proceeds as before, with inequality~\eqref{eq:boundKL} replaced by
\[
|\mu(f_{n+1}-f_n)|\leq \theta^{m}\|f_n\|_\theta+ c\kappa^n\|h^{-1}\|_\infty\|h\|_\infty\|f\|_\cL \left(2\mu(D_p(q^m))\right)^{\frac12}.
\]
We now wish to replace
$\mu(D_p(q^m))$ by 
$\mu(I_m+\log\mu^{\xi_1^{m-1}}_x([x]_\xi))$.
By \cite{sarig2009lecture}[Proposition 3.4], the function $\log h$ satisfies $\var_{n}(\log h)\leq \sum_{\ell>n} \var_\ell(\phi)$.
Equation~\ref{eq:cohom} implies 
\[\var_n{(I_m)}\leq 2\sum_{k\geq n}\var_k(\phi)=:C_n,\] 
which is finite by H\"older summability.
Therefore, as argued in equation~\eqref{eq:approx_m}
\[
|\sum_{s\in S}I_m(sTx)q^m_s(x)-\int I_m d\mu^{\xi_1^{m-1}}_x|\leq C_m\theta^m.
\]
By the properties of the conditional measure, we conclude
\[\mu(D_p(q^m))\leq C_m\theta^m+
\mu(I_m+\log\mu^{\xi_1^{m-1}}_x([x]_\xi)) = C_m\theta^m + P_{GS}(\phi)-(\mu(\phi)+H_\mu(\xi|\xi_1^{m-1}))
\]
where we used the assumptions that $-\mu(\phi)<\infty$ and $H_\mu(\xi|\xi_1^{m-1})<H_\mu(\xi)<\infty$. 

If $\var_\ell(\phi)=0$, $\ell = m$ then in fact $\mu(D_p(q^m))= \mu(I_m+\log\mu^{\xi_1^{m-1}}_x([x]_\xi))$, and by positivity of the Kullback-Leibler convergence,
$P_{GS}(\phi) \geq  \mu(\phi + H_\mu(\xi|\xi_1^{m-1})$ for $m\geq\ell$.
Lemma 2.6 of \cite{kadyrov2016effective} implies that
\[
P_{GS}(\phi)- \left(\mu(\phi)+H_\mu(\xi|\xi_1^{{\lfloor \frac{m-\ell}{2} \rfloor} }))\right)\leq 2 (P_{GS}(\phi)-\frac{1}{m-\ell}\sum_{k=\ell}^{m-1} H_\mu(\xi|\xi_1^{k})).
\]
By additivity of the entropy function
\[
\frac{1}{m-\ell}\sum_{k=\ell}^{m-1} H_\mu(\xi|\xi_1^{k})=\frac{H_\mu(\xi_0^{m-1})-H_\mu(\xi_0^{\ell-1})}{m-\ell} \geq \frac{1}{m} H_\mu(\xi_0^{m-1}) - \frac{\ell}{m-\ell} H_\mu(\xi)  
\]
so that
\[
P_{GS}(\phi)-\left(\mu(\phi)+H_\mu(\xi|\xi_1^{\lfloor\frac{m-\ell}{2}\rfloor})\right)\leq  2 \left( P_{GS}(\phi) - \left( \mu(\phi)+\frac{1}{m} H_\mu(\xi_0^{m-1})  \right) + \frac{\ell}{m-\ell}H_\mu(\xi)  \right) 
\]
Furthermore, $\|f_n\|_\theta=\|f_n-m(f)\|_\theta$, so $\sum_n \|f_n\|_\theta$ is again bounded by $\frac{c}{1-\kappa}\|h^{-1}\|_\infty\|h\|_\infty\|f\|_\cL$. We arrive at Theorem~\ref{thm:finitary_theorem} by choosing $b=(\frac{1}{\sqrt{2}}+\sqrt{2})a$ where $a$ is as in Theorem~\ref{thm:main_theorem}.

\section{Corollaries}
We end this note by stating some corollaries.
\subsection{Approximation by subsystems}
We call a sequence of $T$-invariant probability measures $\{\mu_n\}$ asymptotically $\phi$-equilibrium if $P_{\mu_n}(\phi)\to P_{GS}(\phi)$. These can be constructed by considering a finite subset $S_n\subset S$, and an irreducible submatrix $A_n$ of $A$ supported on $S_n$ so $X_n=\{ x=(x_0,x_1,\dots)\in S_n^{\N_0} : t_{x_ix_j}=1\}$ is $T$-invariant. If $\bigcup_{n}S_n=S$ and $A_n$ is a submatrix of $A_{n+1}$ we call such a sequence exhaustive. Such a sequence always exists (Lemma 3.10 \cite{gurevich1998thermodynamic}). Let us further assume that $A_n$ is $T$-aperiodic and denote $\phi|_n$ the restriction of $\phi$ to $X_n$. Then Ruelle's Perron-Frobenius theorem provides us with a unique $\phi|_n$ equilibrium measure $m_n$ with pressure $P(\phi|_n)=h_{m_n}(T)+m_n(\phi|_n)$. Theorem~\ref{thm:main_theorem} hence tells us that $m_n$ is close to $m$ with bound $\left(P_{GS}(\phi)-P(\phi|_n)\right)^{\frac12}$.
This is a quantification of statement 2) of Theorem 6.3 \cite{gurevich1998thermodynamic} in the setting of Gibbs measures. Indeed, we have $P_{GS}(\phi)-P(\phi|_n)\to0$ which follows from Theorem~2 \cite{sarig1999thermodynamic}
(which in turn generalizes the Markovian case, Lemma 3.10, 2) \cite{gurevich1998thermodynamic}):
\begin{corollary} Let $m$ be a Gibbs measure of a topological mixing countable Markov shift $(X,T)$ for a locally H\"older continuous potential $\phi$. For an exhaustive sequence of aperiodic finite submatrices $A_n$ of $A$ with equilibrium measures $m_n$ we have for any $f\in\cL$
	\[
	|m(f)-m_n(f)|\leq a\|f\|_{\cL}\left(P_{GS}(\phi)-P(\phi|_n)\right)^{\frac12}.
	\]
\end{corollary}

\subsection{Approximation by discrete measures}
In Chapter 7 of \cite{gurevich1998thermodynamic} discrete measures supported on fixpoints of $(T,X)$ are studied. We will make a related analysis, but restrict to fixpoints of the measures $m_n$ just constructed.
Namely, we wish to apply Theorem~\ref{thm:finitary_theorem} to discrete measures supported on fixpoints of $T$ on $X_n$ (which in the case $\phi|_n=0$ is Theorem~1.3 \cite{kadyrov2016effective}):
Let $\text{Fix}_k^n$ denote the fixpoints of order $k$ of $(T,X_n)$. Define $\nu_n^k$ to be the $T$-invariant probability measure on $X_n$ by $\nu_n^k(g)=c_n^k\sum_{y\in \text{Fix}_k^n}e^{(\phi|_n)_k(y)}g(y)$, $c_n^k=\left(\sum_{y\in \text{Fix}_k^n}e^{(\phi|_n)_k(y)}\right)^{-1}$. Let $\xi_n$ be partition of $X_n$ in principal cylinders, then each element in $(\xi_n)_0^{k-1}$ contains at most one fixpoint of order $k$ so that a calculation reveals that
\[
\frac{1}{k}H_{\nu_n^k}((\xi_n)_0^{k-1})+\nu_n^k(\phi|_n)=-\frac{1}{k}\log c_n^k.
\]
In the Markovian case, $\phi(x)=\phi(x_0,x_1)$ the term $(c^k_n)^{-1}$ equals $e^{kP(\phi|_n)}+\sum \lambda^k$ where the sum runs over all eigenvalues $\lambda$ but the top one (equal to $e^{P(\phi|_n)}$) of the matrix $(A_n)_{ij}e^{\phi(i,j)}$ (see Chapter 5 \cite{parry1990zeta}), hence $P(\phi|_n)+\frac{1}{k}\log c_n^k\leq |S_n|\frac{\delta_n^k}{k/2}$ where $\delta_n=\max_{\lambda\neq e^{P(\phi|_n)}} \frac{|\lambda|}{e^{P(\phi|_n)}}<1$ for $k$ (given $n$) large enough.
 Here we are using a first order Taylor approximation of the logarithm.
 Hence $\nu^k_n$ is $b_n(\theta^{\lfloor\frac{k-2}{2}\rfloor}+\left(2|S_n|\frac{\delta_n^k}{k} + \frac{2}{k-2}H_{\nu^k_n}(\xi)\right)^{\frac12})$ close to $m_n$ where $b_n$ is the constant $b$ in Theorem~\ref{thm:finitary_theorem} that depends on the spectral gap and eigenfunction associated to $A_n,\phi|_n$. We summarize:

\begin{corollary} Let $m$ be a Gibbs measure of a topological mixing countable Markov shift $(X,T)$ for a Markovian potential $\phi$ with $\var_1(\phi)<\infty$. For an exhaustive sequence of aperiodic finite submatrices $A_n$  of $A$ and associated state partition $\xi_n$ with equilibrium measures $m_n$ and associated discrete measures $\nu_n^k$ supported on $\text{Fix}_k^n$, we have for all $n$, for all $k\gg_n1$, for any $f\in\cL$
	\[
	|m(f)-\nu_n^k(f)|\leq \|f\|_{\cL}\left(a\left(P_{GS}(\phi)-P(\phi|_n)\right)^{\frac12} + b_n\left(\theta^{k/2}+\left(2|S_n|\frac1k\delta_n^k +\frac{2}{k-2}H_{\nu^k_n}(\xi_n)\right)^{\frac12}\right) \right).
	\]
\end{corollary}

\subsection{Stability of equilibrium measures}
Assume that $(T,X)$ is a topological mixing countable Markov shift with BIP. Let $\phi,\varphi$ be two H\"older continuous potentials with zero Gurevich-Sarig pressure, $P_{GS}(\phi)=P_{GS}(\varphi)=0$. Assume further that the associated Gibbs measures $\mu_\phi,\mu_{\varphi}$ are equilibrium measures, i.e.\ $0=h_{\mu_\phi}(T)+\mu_\phi(\phi)=h_{\mu_\varphi}(T)+\mu_\varphi(\varphi)$. Then by Theorem~\ref{thm:main_theorem}
\[
|\mu_\phi(f)-\mu_{\varphi}(f)|\ll \|f\|_\cL \left(-h_{\mu_\varphi}(T)-\mu_\varphi(\phi)\right)^\frac12=\|f\|_\cL \mu_\varphi(\varphi-\phi)^\frac12\leq \|f\|_\cL\|\phi-\varphi\|_\infty^\frac12.
\]

\begin{corollary} Let $\mu_\phi,\mu_{\varphi}$ be two equilibrium measures of a topological mixing countable Markov shift $(T,X)$ with BIP of two locally H\"older continuous potentials $\phi,\varphi$ with zero pressure. Then for any $f\in\cL$
	\[
	|\mu_\phi(f)-\mu_{\varphi}(f)|\ll \|f\|_\cL\|\phi-\varphi\|_\infty^\frac12.
	\]
\end{corollary}

\bibliographystyle{alpha}
% \bibliography{../bib_rene}

\end{document}